\documentclass[11pt]{article}
\usepackage[utf8x]{inputenc}
\usepackage{amsfonts,amssymb,amsmath,amscd}
\usepackage{enumitem}
\usepackage{upgreek}

\usepackage[alphabetic,msc-links]{amsrefs}
\usepackage{fullpage}

\usepackage{hyperref} 
\hypersetup{pdfborder = 0 0 0}

\newcommand{\affil}{\thanks}

\setlist[itemize]{noitemsep, topsep=0pt}

\usepackage{xcolor}
\definecolor{keywordcolor}{rgb}{0, 0.1, 0.6}   
\definecolor{commentcolor}{rgb}{0.4, 0.4, 0.4}   
\definecolor{symbolcolor}{rgb}{0, 0.0, 0.0}    
\definecolor{tacticcolor}{rgb}{0, 0.3, 0.1}    
\definecolor{sortcolor}{rgb}{0, 0, 0}      
\definecolor{errorcolor}{rgb}{1, 0, 0}           
\definecolor{stringcolor}{rgb}{0.5, 0.3, 0.2}    

\usepackage{listings}

\usepackage{upgreek}

\def\lean{\lstinline[language=lean]}

\usepackage{tikz}
\usetikzlibrary{tikzmark}

\title{
	Anatomy of a Formal Proof
}

\author{
  Jeremy Avigad
  \affil{
	  Carnegie Mellon, USA. Professor. \texttt{avigad@cmu.edu}
    }
  \and
  Johan Commelin
  \affil{
	  Utrecht University, The Netherlands. Assistant Professor.
	  \texttt{j.m.commelin@uu.nl}
    }
  \and
  Heather Macbeth
  \affil{
	  Fordham University, USA. Assistant Professor.
	  \texttt{hmacbeth1@fordham.edu}
    }
  \and
  Adam Topaz
  \affil{
	  University of Alberta, Canada. Associate Professor.
	  \texttt{topaz@ualberta.ca}
    }
}

\begin{document}

\maketitle

\begin{abstract}
Interactive proof assistants make it possible for ordinary mathematicians to write definitions and theorems in a formal proof language, like a programming language, so that a computer can parse them and check them against the rules of a formal axiomatic foundation. This article describes the experience of working with a proof assistant and considers the impact the technology will have on mathematics.
\end{abstract}

\section*{Introduction}

It has been a long day and you are making your way through a paper related to your work. You suddenly come across the following remark: ``\ldots since $x$ and $y$ are eigenvectors of $f$ with distinct eigenvalues, they are linearly independent.'' Wait---how does the proof go? You should really know this. Here $x$ and $y$ are nonzero elements of a vector space $V$ and $f : V \to V$ is a linear map. You force yourself to pick up a pen and write down the following argument:
\begin{quote}
Let $f(x) = \mu x$ and $f(y) = \nu y$ with $\mu \ne \nu$. Suppose $a x + b y = 0$. Applying $f$ and using linearity, we have $\mu a x + \nu b y = 0$. Multiplying the original equation by $\nu$, we have $\nu a x + \nu b y = 0$. Subtracting the two yields $(\mu - \nu) a x = 0$ and since $\mu - \nu$ and $x$ are nonzero, we have $a = 0$. The corresponding argument with $x$ and $y$ swapped yields $b = 0$, so the only linear combination of $x$ and $y$ that yields $0$ is the trivial one.
\end{quote}
That works, doesn't it?

Your colleagues have all gone home and there is nobody around to discuss this with. So, instead, you turn to your computer and start up Lean, the proof assistant you happen to use. Can you prove the claim formally? As you type, the information window in your editor complains about every misstep---the syntax is fiddly, and you have to get the notation and the instructions just right---but that's okay, because working through the proof is relaxing and kind of fun. Lean often makes you spell out arguments that are painfully obvious, but you have found that if you set things up just right, it will cheerfully fill in some of the details. After a short while, you have success! Lean signs off on the proof, confirming you have managed to construct a formal derivation in the system's axiomatic foundation.

\begin{lstlisting}[numbers=right]
import Mathlib.LinearAlgebra.LinearIndependent

variable [Field K] [AddCommGroup V] [Module K V]

example (f : V →ₗ[K] V)
    (μ ν : K) (hμν : μ ≠ ν)
    (x y : V) (hx₀ : x ≠ 0) (hy₀ : y ≠ 0)
    (hx : f x = μ • x) (hy : f y = ν • y) :
    ∀ a b : K,
      a • x + b • y = 0 → a = 0 ∧ b = 0 := by
  intro a b hab
  have :=
  calc (μ - ν) • a • x
      = (a • μ • x + b • ν • y) -
        ν • (a • x + b • y) := by module
    _ = f (a • x + b • y) -
        ν • (a • x + b • y) := by simp [hx, hy]
    _ = 0 := by simp [hab]
  simp_all [sub_eq_zero]
\end{lstlisting}
What a lovely proof! It took some effort to work it out, but it was a pleasure seeing the steps play out formally, and you are proud of the result.

\section*{Understanding the proof}

Let's take a look at what you have done. You are using Lean's mathematical library, Mathlib, a communally developed and maintained repository of formally verified mathematics.
At the time of writing, Mathlib contains more than 80,000 definitions and more than 160,000 theorems, starting from axiomatic primitives and developing mathematics from the bottom up.
It also contains dozens of user-contributed automated reasoning procedures that help with the formalization process, as well as notation declarations and configuration information, all of which codify different aspects of our mathematical understanding.

To avoid having to load the entire library into memory at once, Lean asks you to tell the system what parts of the library you want to use.
So the first line of the proof, the one that starts with the keyword \lstinline{import}, tells Lean you want to use the file with the name shown.
That pulls in all the files in the library that those notions depend on, which is quite a lot, including basic algebra, properties of scalars and vectors, and so on.

The next line, the one that begins with the word \lstinline{variable}, declares some of the objects you want to work with: a field, $K$, and a vector space, $V$.
You could have, alternatively, put these declarations after the \lstinline{example} keyword
on line 5; using the \lstinline{variable} command to declare them separately is especially
convenient when multiple definitions and theorems share the same data and hypotheses.
It's a quirk of the library's design that you declare that $V$ is an additive, commutative, group, coupled with a scalar action that turns $V$ into a $K$-module.
(Remember, a vector space is nothing more than a module over a ring in which the ring in question is actually a field.)
Mathlib includes thousands of axiomatically declared structures, and the contributors to the library take great care to maximize reusability by breaking structural hypotheses into reusable pieces that can easily be configured and composed.

Next comes the keyword \lstinline{example} (line~5), which indicates that the result you are proving is not intended for future use.
Using \lstinline{example} is good for experimentation.
To prove a lemma or theorem that you intend to use later, you would instead use \lstinline{lemma} or \lstinline{theorem} and provide a name.
In that case, once Lean has processed the proof, it is stored in the \emph{environment}, which means that any file that imports this one can see them.
If you think you have proved a theorem that ought to go into Mathlib, you can issue a \emph{pull request} to add it to the library.
Mathlib is overseen by teams of \emph{maintainers} and \emph{reviewers} who moderate and update the contents, ensuring the quality and stability of the library.

You begin stating the claim by introducing the key players:
\begin{itemize}
  \item a linear map $f$ from $V$ to itself. Lean's notation for such linear maps is \lstinline{V →ₗ[K] V},
  \item two elements of $K$, denoted $\mu$ and $\nu$, which will end up being your eigenvalues, and
    a hypothesis called \lstinline{hμν} that says $\mu$ and $\nu$ are distinct,
  \item two vectors called $x$ and $y$, as well as hypotheses \lstinline{hx₀} and \lstinline{hy₀} that say that they are nonzero,
  \item and, finally, your two key assumptions, \lstinline{hx} and \lstinline{hy}, which say that $f(x) = \mu \cdot x$ and $f(y) = \nu \cdot y$.
\end{itemize}
In the notation \lstinline{V →ₗ[K] V}, the \lstinline{ₗ} is a fixed symbol signifying linear maps,
whereas the parameter in brackets indicates that you mean linearity over \lstinline{K}.

At this point, you have set up your \emph{context}, that is, the data and hypotheses you need to state your claim.
You have your $K$-vector space $V$, a linear map $f : V \to V$, two scalars $\mu$ and $\nu$, and eigenvectors $x$ and $y$ for $\mu$ and $\nu$ respectively.
At that point you are ready to state the conclusion (lines~9 and~10),
namely, that $x$ and $y$ are linearly independent: \lstinline{∀ a b : K, a • x + b • y = 0 → a = 0 ∧ b = 0}.
You have decided to use the definition of linear independence for a pair of vectors: Whenever $a$ and $b$ are scalars, if $a \cdot x + b \cdot y = 0$ then $a = 0$ and $b = 0$.

Notice that you have used the symbol \lstinline{:} in expressions like \lstinline{x y : V} (line~7) to express that \lstinline{x} and \lstinline{y} are elements of \lstinline{V},
rather than writing \lstinline{x ∈ V} and \lstinline{y ∈ V}.
In Lean's axiomatic foundation, a version of \emph{dependent type theory}, every object has a fundamental data type.
In this example, we would say that that $x$ and $y$ are \emph{terms} of the \emph{type} $V$.
Interestingly, we use the same notation in an expression like \lstinline{hx₀ : x ≠ 0} (also line~7) to express that \lstinline{hx₀} is a label for the assumption \lstinline{x ≠ 0}.
In dependent type theory, assertions like \lstinline{x ≠ 0} (also known as \emph{propositions}) are analogous to data types, and they are handled by the same fundamental mechanisms.
Lean checks that an expression like \lstinline{a • x + b • y} is a well-formed expression of type \lstinline{V}, given information about the variables and symbols involved, and, in the same way, it checks that a formal expression is a well-formed proof of the proposition \lstinline{∀ a b : K, a • x + b • y = 0 → a = 0 ∧ b = 0}, given the data and assumptions it depends on.
That is why the latter appears after the colon in the example's conclusion: the statement of the example announces your intent to construct a formal proof of the conclusion, given the data and assumptions that come before.

What comes after the symbol \lstinline{:=} (line~10) is the proof itself, or, more precisely, instructions that tell Lean how to construct the proof.
If you are successful, the corresponding expression is stored in memory and checked for correctness by Lean's trusted \emph{kernel}.
The keyword \lstinline{by} after the \lstinline{:=} instructs Lean to enter \emph{tactic mode}, which means that the text that follows should be interpreted as a list of instructions that tell Lean how to build the required proof.
A \emph{tactic} is a procedure that automatically fills in a chain of formal inferences that
is needed to justify a reasoning step.
Such a chain can be as short as a single logical axiom or rule, but it
can also be quite long and involved.
Tactics thus play an important role in bridging the gap between the kind of reasoning that
is intuitively clear and natural to mathematicians and the stringent axiomatic rules
embodied by a formal system.

The first tactic you use is \lstinline{intro a b hab} (line~11).
This introduces the two variables \lstinline{a : K}, \lstinline{b : K} which your statement quantifies over, as well as the antecedent of the implication, that is,
the assumption \lstinline{hab : a • x + b • y = 0}.
As you type or move your cursor around a proof, Lean displays the \emph{tactic state}, i.e.~information that is relevant at that point in the proof, in its \emph{infoview} window.
After the \lstinline{intro} tactic, the tactic state looks something like this:

\begin{lstlisting}
K : Type
V : Type
f : V →ₗ[K] V
μ ν : K
hμν : μ ≠ ν
x y : V
hx₀ : x ≠ 0
hy₀ : y ≠ 0
hx : f x = μ • x
hy : f y = ν • y
a b : K
hab : a • x + b • y = 0
⊢ a = 0 ∧ b = 0
\end{lstlisting}
This is nice a summary of where you are in the proof, including the objects and assumptions you started with as well as the objects \lstinline{a}, \lstinline{b} and the assumption \lstinline{hab} that you introduced in the first step.
The line that starts with \lstinline{⊢} indicates that your current \emph{goal} is to prove the conjunction \lstinline{a = 0 ∧ b = 0}.

While writing the proof, you notice that Lean complains with warnings and error messages.
This is expected, since the proof is incomplete.
You can appease Lean by apologizing for the incomplete proof:
If you use the \lean"sorry" tactic on the final line of the proof,
Lean will not raise an error over the fact that the proof is incomplete,
but it will still raise a gentle warning about the use of \lean"sorry".

The next step in the proof starts with \lean"have := calc" (line~12).
This introduces a calculation, similar to a \texttt{\textbackslash{}begin\{align*\} ... \textbackslash{}end\{align*\}} block in \LaTeX.
The calculation itself is very similar to the one in the proof sketch above,
but notice that the individual steps in the calculation are justified by short subproofs like \lean"by module" (line~15) or \lean"by simp [hx, hy]" (line~17).
The \lean"module" tactic proves equalities of universal linear expressions,
in other words, linear equalities that are true in all modules, and do not use specific facts about the module or ring at hand.
The \lean"simp" tactic is a powerful tool that uses a database of equalities and equivalences from Mathlib,
together with user-specified equations (like \lean"hx, hy") to rewrite its goal into a simpler form.
Happily, in this case, the goal becomes a trivial equality, and the goal is therefore closed.
In a moment, we will consider the steps that \lean"simp" has taken in greater detail.

After this calculation,
you completed your proof with the \lean"simp_all" tactic (line~19),
which is a variant of \lean"simp" that recursively uses all the hypotheses to simplify all hypotheses and the goal.
Now if you leave a \lean"sorry" after the \lean"simp_all",
Lean will complain that there is nothing to be sorry about
and insist that you remove the apology.
Lean checks the correctness of the proof and confirms it with its silence:
the absence of errors is the proclamation that you have succeeded.

\section*{Using automation}

So what exactly does the \lean"simp" tactic do?
The following provides a manual, more detailed proof of one step of the calculation in your proof, with the direction of the equality reversed.
\begin{lstlisting}
variable (f : V →ₗ[K] V)
  (μ ν : K) (x y : V) (a b : K)

example (hx : f x = μ • x) (hy : f y = ν • y) :
  f (a • x + b • y) = (a • μ • x + b • ν • y) := by
calc f (a • x + b • y)
    = f (a • x) + f (b • y) := by
        rw [map_add]
  _ = a • f x + b • f y := by
        rw [map_smul, map_smul]
  _ = (a • μ • x + b • ν • y) := by
        rw [hx, hy]
\end{lstlisting}
This example demonstrates the \emph{rewrite} tactic, denoted \lean"rw".
In the calculation, we first rewrite with the lemma \lean"map_add",
which states that for any linear map \lean"f",
the equation \lean"f (v + w) = f v + f w" holds.
The next step in the calculation is justified
by rewriting twice with the lemma \lean"map_smul",
which states that for any linear map \lean"f", scalar \lean"c", and vector \lean"v",
the equation \lean"f (c • v) = c • f v" holds.
(Here \lean"smul" stands for ``scalar multiplication.'')
The underscores before each subsequent step are part of the \lstinline{calc}
syntax, indicating to the parser that the calculation continues.
The proof concludes by rewriting with the hypotheses \lean"hx" and \lean"hy",
which assert that \lean"x" and \lean"y" are eigenvectors of \lean"f" with eigenvalues \lean"μ" and \lean"ν" respectively.

The lemmas \lean"map_add" and \lean"map_smul" are part of the Mathlib library,
and they are labeled with the \emph{attribute} \lean{@[simp]}.
This attribute tells Lean that the lemma should be added to the database of lemmas that \lean"simp" uses to simplify expressions.
And that is the reason why the \lean"simp" tactic could prove the goal in one go:
it chained together the \lean"rw" steps that we spelled out step by step in the calculation above.

One of the challenges of formalization is that we often need to spell out, in painful detail, inferences that seem obvious or straightforward to us.
The more we can get the computer to fill in, the better.
The simplifier is an important example of \emph{automated reasoning} that can help in this respect.
The \lstinline{module} tactic is another.

Broadly speaking, there are two classes of automation.
Firstly, there is \emph{general purpose} automation.
One example that we have seen is the \lean"simp" tactic, the simplifier.
Another example is \lean"aesop",
a tactic that provides ``Automated Extensible Search for Obvious Proofs.''
In the proof assistant Isabelle, there is a tool called \lean"sledgehammer",
that can search for proofs using a large database of lemmas.
And recently AI copilots have demonstrated the ability to suggest tactics and fill in parts of proofs in Lean.
The wide applicability of these general-purpose tools is balanced by the fact that, at least for now, they can only assist with proofs that are relatively straightforward.

The other class of automation is \emph{domain specific} automation.
Here we have also seen an example: the \lean"module" tactic.
Other examples include the \lean"ring" tactic, which solves equations in commutative rings,
the \lean"linarith" tactic, which solves linear arithmetic problems,
and the \lean"fun_prop" tactic, which proves that functions satisfy a given property such as ``continuity'' or ``measurability.''
These tactics are less general; they typically have a well-defined and narrow scope,
but they can be very powerful in their domain of applicability.
A prime example of this approach comes from the use of software tools known as
SAT solvers and SMT solvers:
if a claim can be encoded as a boolean formula or a formula in some decidable theory,
then these can be used to justify it automatically.
This has been fruitfully applied in non-trivial ways, see for example~\cite{MR4449705}.

\section*{Generalizing the hypotheses}


Often, one of the first things we do when we have proved a lemma or a theorem is check whether the hypotheses can be weakened, in order to increase its applicability.
Proof assistants are especially helpful in this respect, because they enable us to tinker with hypotheses interactively and see what breaks.
Returning to our example, upon verifying the initial result, you might wonder:
what are the minimal assumptions that we need to make this proof go through?
Let's experiment.
The proof never mentions inverses of scalars,
so it should work for a larger class of rings.
In a first attempt, we might try to replace the vector space \lstinline{V}
by an arbitrary module \lstinline{M} as follows:
\begin{lstlisting}
variable {R M : Type}
  [Ring R] [AddCommGroup M] [Module R M]
\end{lstlisting}
But Lean complains at
the first step in the calculation, which requires
\lean"b • ν • y = ν • b • y".
We can address this by assuming
that the ring of scalars is commutative.
However, after making that change, the final step of the proof is still broken.
A bit of reflection shows that this step uses
the result of the calculation steps together with the assumptions
\lean"hμν : μ ≠ ν" and \lean"hx₀ : x ≠ 0"
to prove that \lean"a = 0".
In other words, we need the additional property that
\lstinline{r • m = 0} implies
that \lstinline{r = 0} or \lstinline{m = 0} for every \lstinline{r} in
\lstinline{R} and \lstinline{m} in \lstinline{M}.
Mathlib expresses this as the property \lean"NoZeroSMulDivisors R M",
and we can add that assumption as follows:
\begin{lstlisting}
variable {R M : Type}
  [CommRing R] [AddCommGroup M] [Module R M]
  [NoZeroSMulDivisors R M]

example (f : M →ₗ[R] M)
  (μ ν : R) (hμν : μ ≠ ν)
  (x y : M) (hx₀ : x ≠ 0) (hy₀ : y ≠ 0)
  (hx : f x = μ • x) (hy : f y = ν • y) :
  ∀ a b : R, a • x + b • y = 0 → a = 0 ∧ b = 0
\end{lstlisting}
Lean accepts this statement with the same proof as before!
We have therefore obtained a more general theorem without changing a single character of the proof.


\section*{Strengthening the theorem}

Let's see if we can strengthen your result even further.
You started with a basic version involving vector spaces over fields and only two eigenvectors,
and we just generalized it to suitable modules over commutative rings.
Where can we go from here?

Let's try to generalize the theorem to
arbitrary families of eigenvectors.
To do this, we will need to think carefully about how to express,
in Lean, the fact that an arbitrary family of vectors is linearly independent.
We could write down a definition from scratch, but it makes sense to see if we can take advantage of things that are already in the library.
If we navigate to the Mathlib documentation webpage%
\kern-.06em\footnote{\url{https://leanprover-community.github.io/mathlib4_docs/}}
and start typing \texttt{linear independent} in the search box,
one of the first few results that comes up is the following:
\begin{lstlisting}
def LinearIndependent {ι : Type}
  (R : Type) {M : Type}
  (v : ι → M) [Semiring R]
  [AddCommMonoid M] [Module R M] :
  Prop
\end{lstlisting}
Here \lean"R" and \lean"M" are the relevant ring and module, as before, \lean"ι" is an indexing type, and \lean"v" is a family of elements of
\lean"M" indexed by \lean"ι".
It is standard, and convenient, to
represent a family $(v_i)_{i \in \iota}$ of elements of $M$ as a function \lstinline{v : ι → M},
in which case the $i$th element, $v_i$, is simply written \lstinline{v i}.
All the arguments in curly and square brackets are generally left implicit, which is to say, we expect to write \lean"LinearIndependent R v" and have Lean figure out the rest.
The annotation \lean"Prop" means that the expression
\lean"LinearIndependent R v" is a proposition, namely, the proposition that the family \lean"v" is linearly independent over \lean"R".
At the moment, we don't need to know the body of the definition; we can use it as a black box.

We will add an indexing type \lean"ι" to the statement of our theorem, and since we also want to say that the
vectors are all \emph{eigenvectors}, we will also use \lean"ι" to index a family $\mu$ of scalars which will act as the eigenvalues.
We can now formulate the statement we are after as follows:
\begin{lstlisting}
example {ι : Type} (f : M →ₗ[R] M)
    (μ : ι → R) (hμ : Function.Injective μ)
    (v : ι → M) (hv : ∀ i, v i ≠ 0)
    (h : ∀ i, f (v i) = μ i • v i) :
    LinearIndependent R v
\end{lstlisting}
Note that we have formulated the fact that we are considering \emph{distinct} eigenvalues by assuming $\mu$ is injective.

We should check to see whether something like this theorem is already
in the library.
If we go back to the documentation webpage and search for \texttt{eigenvector linear independent},
only two results come up:
\begin{lstlisting}
theorem
  Module.End.eigenvectors_linearIndependent
  (f : Module.End R M)
  (μs : Set R) (xs : ↑μs → M)
  (h_eigenvec : ∀ (μ : ↑μs),
    f.HasEigenvector (↑μ) (xs μ)) :
  LinearIndependent (ι := ↑μs) R xs
\end{lstlisting}
and the variant
\begin{lstlisting}
theorem
  Module.End.eigenvectors_linearIndependent'
  {ι : Type} (f : Module.End R M)
  (μ : ι → R) (hμ : Function.Injective μ)
  (v : ι → M)
  (h_eigenvec : ∀ (i : ι),
    f.HasEigenvector (μ i) (v i)) :
  LinearIndependent R v
\end{lstlisting}
The difference between the two is that the first is about a set of scalars and a function assigning an eigenvector to each scalar, whereas the second one is about indexed families of scalars and eigenvectors.
The second one is more promising for our application, since we also chose to use indexed families of vectors and scalars.
To use this theorem to prove our version, we invoke the
\lean"apply" tactic.
\begin{lstlisting}
example {ι : Type} (f : M →ₗ[R] M)
    (μ : ι → R) (hμ : Function.Injective μ)
    (v : ι → M) (hv : ∀ i, v i ≠ 0)
    (h : ∀ i, f (v i) = μ i • v i) :
    LinearIndependent R v := by
  apply Module.End.eigenvectors_linearIndependent'
\end{lstlisting}
This leaves us with a number of goals:
\begin{enumerate}
  \item \lean"Function.Injective ?μ", asking us to prove that something is injective.
  \item \lean"∀ (i : ι), ?f.HasEigenvector (?μ i) (v i)", which looks like it should have something to do with our assumptions \texttt{h} and \texttt{hv}.
  \item \lean"Module.End R M", asking us to provide an endomorphism of $M$; this should just be our linear map \texttt{f}.
  \item And finally, \lean"ι → R", which will be our $\mu$.
\end{enumerate}
The \lean"?μ" in the first two goals and the \lean"?f" in the second goal mean that Lean does not yet know how to instantiate the variables $\mu$ and $f$ in the theorem we have invoked.
It left those tasks as the third and fourth goals; Lean expects that it will be more convenient for us to provide that information implicitly when we solve the other goals.

We have made progress, but we still need to provide the information requested.
We should clearly use \lean"hμ" to solve the first goal, which we do by adding the next line to our proof:
\begin{lstlisting}
  apply Module.End.eigenvectors_linearIndependent'
  exact hμ
\end{lstlisting}
As you might guess, the \lean"exact" tactic tells Lean to use the assumption \lean"hμ" to close the goal.
As a side effect, that also closes goal 4: Lean is now able to infer that the family in question is $\mu$.
We then have two goals left:
\begin{enumerate}
  \item \lean"∀ (i : ι), ?f.HasEigenvector (μ i) (v i)", and
  \item \lean"Module.End R M".
\end{enumerate}
We know that the endomorphism of $M$ should be $f$, so we swap the order of the goals and give Lean
this information:
\begin{lstlisting}
  apply Module.End.eigenvectors_linearIndependent'
  exact hμ
  swap ; exact f
\end{lstlisting}
This leaves us with one last goal, namely, \lean"∀ (i : ι), Module.End.HasEigenvector f (μ i) (v i)".
We need to use our assumptions \texttt{h} and \texttt{hv}, but it's unclear how to package them together to satisfy the definition of \lean"Module.End.HasEigenvector".
Instead of going back to the documentation page and looking up the definition, we can ask Lean what we have to do.
\begin{lstlisting}
  apply Module.End.eigenvectors_linearIndependent'
  exact hμ
  swap ; exact f
  intro i ; constructor
\end{lstlisting}
As before, the \lean"intro" tactic introduces an arbitrary \lean"i",
and then the \lean"constructor" tactic tells Lean that we are ready to
provide the information needed to show that \lean"v i" is an eigenvector of $f$ with
eigenvalue \lean"μ i".
This makes progress, but we still have the following goals to fulfill:
\begin{enumerate}
  \item \lean"v i ∈ Module.End.eigenspace f (μ i)"
  \item \lean"v i ≠ 0", which is obviously an application of \texttt{hv}.
\end{enumerate}
For 1, the hypothesis \lean"h : ∀ i, f (v i) = μ i • v i" should do the trick, but applying it directly doesn't work.
After digging into the details, we see that this is because the eigenspace of $f$ with respect to a scalar $\mu$ is defined as the kernel of $f - \mu \cdot 1$, so we will need to convert \texttt{h} to this form.
We could do this manually, but we can first check whether we can use a preexisting lemma from the library.
The fastest way to do this, in the middle of a Lean proof, is to try the \lean"exact?" tactic.
This will do a search for ways of applying existing lemmas to close the goal \emph{exactly}.
This tactic doesn't always work, but it can't hurt to try.
Aha! In this case it tells us right away that we can close the first goal by writing \lean"exact Module.End.mem_eigenspace_iff.mpr (h i)".
The second goal is also easy to close by applying \texttt{hv}, so now we have a complete proof:
\begin{lstlisting}
  apply Module.End.eigenvectors_linearIndependent'
  exact hμ
  swap ; exact f
  intro i
  constructor
  exact Module.End.mem_eigenspace_iff.mpr (h i)
  apply hv
\end{lstlisting}
The beginning of our proof was a bit messy, but we can clean things up to obtain a nice final result:
\begin{lstlisting}
example {ι : Type} (f : M →ₗ[R] M)
    (μ : ι → R) (hμ : Function.Injective μ)
    (v : ι → M) (hv : ∀ i, v i ≠ 0)
    (h : ∀ i, f (v i) = μ i • v i) :
    LinearIndependent R v := by
  apply
    Module.End.eigenvectors_linearIndependent' f μ hμ
  intro i
  constructor
  · exact Module.End.mem_eigenspace_iff.mpr (h i)
  · apply hv
\end{lstlisting}
The \lean"apply" tactic now supplies \lean"f", \lean"μ", and \lean"hμ" right away, leaving
only one remaining goal.

We have thus achieved the level of generality we were after.
Should this theorem be added to Mathlib?
Upon consideration, we should conclude that our theorem is not substantially different from \lean"Module.End.eigenvectors_linearIndependent'". In our proof, we only added some plumbing and changed the way we talk about eigenvectors to match
the theorem in the library.
Now that we know about \lean"Module.End.HasEigenvector", it seems that this is the standard way to talk about eigenvectors in Mathlib.
If we were to modify our statement to use this instead of our bespoke \texttt{h} and \texttt{hv}, our theorem would be a direct application of \lean"Module.End.eigenvectors_linearIndependent'".
In other words, we have come to realize that the theorem in the library
is what we really wanted all along.
But the good news is that we have learned a lot in the process,
and we are now much more comfortable reasoning about linear independence and eigenvectors with Mathlib.







\section*{Using version control}


Programmers and computer scientists have long made use of
\emph{version control} platforms like GitHub to work on large, collaborative software projects.
Formalization has brought the same tools and methodologies to mathematics.
At the time of writing, Mathlib comprises approximately 5,000 files and 1.5 million lines,
written by over 300 contributors.
These contributions go through an open review process on GitHub,
before they are merged into the main repository.
In total there have been a bit more than 30,000 contributions
since Mathlib's inception in 2017.

A nice thing about version control tools is that they maintain the entire history of the project, allowing us to see what has changed and when.
For example,
we can trace the history of the concept of ``linear independence'' in Mathlib.
It all started on December 7, 2017,
when Johannes H\"olzl committed%
\footnote{\url{https://github.com/leanprover-community/mathlib3/commit/c32d01d}}
the file \lean"algebra/linear_algebra/basic.lean"
to the repo\-si\-tory.
The file was 708 lines long, and on line 186 it contained a definition of linear independence.
\begin{lstlisting}
def linear_independent (s : set β) : Prop :=
  ∀ l : lc α β, (∀x∉s, l x = 0) →
    l.sum (λv c, c • v) = 0 → l = 0
\end{lstlisting}
The terms \lean"l" quantify over \lean"lc α β",
the type of all linear combinations of elements of \lean"β" with coefficients in \lean"α".

On March 10, 2018, an administrative operation%
\footnote{\url{https://github.com/leanprover-community/mathlib3/commit/d010717}}
moved \lean"linear_algebra/" out of \lean"algebra/" so that it became a top-level folder.
Another such move occurred on January 15, 2019,
when Simon Hudon moved%
\footnote{\url{https://github.com/leanprover-community/mathlib3/commit/78f1949}}
all the mathematical content into \lean"src/"
to separate them from the tests and other auxiliary files.

On July 3, 2019, Alexander Bentkamp morphed%
\footnote{\url{https://github.com/leanprover-community/mathlib3/commit/d2c5309}}
the definition into
\begin{lstlisting}
def linear_independent : Prop :=
  (finsupp.total ι M R v).ker = ⊥
\end{lstlisting}
In other words, a collection of vectors \lean"v" in \lean"M" is linearly independent
if the natural map from the free module generated by the vectors \lean"v" to \lean"M" has trivial kernel.

Then, on October 5, 2020,
Anne Baanen split%
\footnote{\url{https://github.com/leanprover-community/mathlib3/commit/2127165}}
linear independence and the accumulated supporting theory into a separate file:
\lean"src/linear_algebra/linear_independent.lean".
The file was 918 lines long.
The story continued on February 23, 2023,
when Pol'tta / Miyahara Kō ported%
\footnote{\url{https://github.com/leanprover-community/mathlib4/commit/db8c8ed}}
the file to Lean 4
as part of a massive collaborative effort to move all of Mathlib to the new version of Lean.
The filename is now \lean"Mathlib/LinearAlgebra/LinearIndependent.lean",
and it lives in the new \lean"mathlib4" repository on GitHub.
The definition has not changed substantially since the change by Bentkamp.
At the time of writing, it reads as follows:
\begin{lstlisting}
def LinearIndependent : Prop :=
  LinearMap.ker (Finsupp.total ι M R v) = ⊥
\end{lstlisting}

\section*{Conclusions}

The earliest programs for checking mathematical proofs include Nicolaas de Bruijn's Automath system, launched in 1967, and Andrzej Trybulec's Mizar system, launched in 1973. Since then, dozens of proof assistants have been developed; Coq, Isabelle, and HOL Light are among the more prominent ones still in use today.
The Lean project, launched by Leonardo de Moura in 2013, is a relative newcomer.
A special issue of the \emph{Notices}, with articles by John Harrison~\cite{MR2463992}, Thomas Hales~\cite{MR2463990}, and Freek Wiedijk~\cite{MR2463993}, surveyed the state of the field in 2008.
We have come a long way since then.

Proof assistants are now commonly used in industry to verify hardware, software, network protocols, cryptographic protocols, cyberphysical systems, and more.
Mathematicians have only recently begun to embrace the technology,
and it is becoming clear that there are several benefits to representing mathematics in digital form.
Just as the word processor opened up new opportunities for written expression and communication, the digitization of mathematics opens up new opportunities for mathematical research and teaching.\footnote{See \cites{massot2021formalize,MR4726989,MR4680264} and the recent pair of special issues for the \emph{Bulletin of the American Mathematical Society} (Volume 61, Numbers 2 and 3), ``Will Machines Change Mathematics?''} The practical benefits are not the only motivation; to many of us, formalizing mathematics feels like the right thing to do. Mathematical definitions and theorems \emph{deserve} to be rendered digitally.

It is also becoming clear that the technology is here to stay. Mathlib currently has roughly 1.5 million lines of code and continues to grow.
Important results, including
foundations for Clausen and Scholze's condensed mathematics,\footnote{\url{https://leanprover-community.github.io/liquid/}}
the polynomial Freiman--Ruzsa conjecture,\footnote{\url{https://teorth.github.io/pfr/}}
and an exponentially improved upper bound to Ramsey's theorem,\footnote{\url{https://xenaproject.wordpress.com/2023/11/04/formalising-modern-research-mathematics-in-real-time/}}
have been formally verified before journal referees signed off on them.
A number of collaborative verification projects have been launched,
including a proof of the sphere eversion theorem\footnote{\url{https://leanprover-community.github.io/sphere-eversion/}}
and a proof of a strengthened version of
Carlson's theorem on pointwise almost everywhere convergence of Fourier series.\footnote{\url{http://florisvandoorn.com/carleson/}}

We expect that, in the years to come, AI copilots that combine neural and symbolic methods will significantly ease the burden of formalization.
More dramatically, we expect that the technology we have discussed here
will play a significant role in the discovery of new mathematics.
Note that DeepMind's AlphaProof, which was deemed to have performed at the level of a silver medalist at the most recent International Mathematical Olympiad,%
\footnote{\url{https://deepmind.google/discover/blog/ai-solves-imo-problems-at-silver-medal-level/}}
was trained to find formal proofs in Lean.
It is exciting to think about what a synergetic combination of machine learning, symbolic methods, and user interaction will bring to mathematics in the years ahead.

You can find online documentation and tutorials for all the proof assistants we have just mentioned.
Lean's lively social media channel on the Zulip platform is welcoming to newcomers,
and the Lean community web pages contain links\footnote{\url{https://leanprover-community.github.io/learn.html}} to learning resources, like the Natural
Number Game, to help you get started.
Proof assistants are not easy to use, and learning to formalize mathematics requires
significant time and effort.
Interaction with proof assistants like Lean comes naturally, however, to those who have grown up immersed in computational technology.
We have therefore found that one of the best ways to take advantage of proof assistants is to have our students help us out.
Formal mathematics is a language, and only they can claim to be among the first generation of native speakers, while the rest of us struggle to master the grammar and intonation.

With all the changes looming, we ought to be concerned about the ways that proof assistants and AI will change the mathematics that we know and love.
It is therefore all the more important for those of us who are more settled in our careers to play an active role in the adoption of the new technologies, using our mathematical values and expertise to guide our students as they negotiate the changing landscape.
The new developments offer us a wonderful opportunity to lead from behind, and it falls on all of us to support the next generation of mathematicians as they forge a path into the digital future.

\section*{Acknowledgments}

We are grateful to the Hausdorff Research Institute for Mathematics for hosting three
of us for the trimester program, ``Prospects of Formal Mathematics,'' in the summer of
2024, during which most of this article was written.
We are also grateful to Paul Buckingham and three anonymous referees for helpful comments,
corrections, and suggestions.

\bibliography{refs}

\end{document}